\documentclass[letterpaper, 10 pt, conference, sort&compress]{ieeeconf}

\IEEEoverridecommandlockouts

\title{\LARGE \bf
	Numerical optimal control for delay differential equations:\\A simultaneous approach based on linearization of the delayed state
}

\author{Tobias K. S. Ritschel and Søren Stange
	\thanks{*T. K. S. Ritschel is with Department of Applied Mathematics and Computer Science, Technical University of Denmark, DK-2800 Kgs. Lyngby, Denmark
		{\tt\small tobk@dtu.dk}}%
}

\usepackage{amsmath}
\usepackage{xcolor}
\usepackage{graphicx}
\usepackage{epstopdf}

\newcommand{\pdiff}[2]{\frac{\partial{#1}}{\partial{#2}}}

\newcommand{\incr}{\,\mathrm{d}}

\begin{document}
	\maketitle
\thispagestyle{empty}
\pagestyle{empty}

	\begin{abstract}
	Time delays are ubiquitous in industry, and they must be accounted for when designing control strategies. However, numerical optimal control (NOC) of delay differential equations (DDEs) is challenging because it requires specialized discretization methods and the time delays may depend on the manipulated inputs or state variables. Therefore, in this work, we propose to linearize the delayed states around the current time. This results in a set of implicit differential equations, and we compare the steady states and the corresponding stability criteria of the DDEs and the approximate system. Furthermore, we propose a simultaneous approach for NOC of DDEs based on the linearization, and we discretize the approximate system using Euler's implicit method. Finally, we present a numerical example involving a molten salt nuclear fission reactor.
\end{abstract}

	\section{INTRODUCTION}\label{sec:introduction}
Time delays are ubiquitous in industrial processes, and they can significantly affect the stability of a system. They arise due to, e.g., transport in pipes (advection), slow mixing (diffusion), growth kinetics, and communication delays, and they must be accounted for when designing control systems. Mathematical models are often used in the development of control strategies, e.g., for analysis, prediction, monitoring, control, and optimization, and the time delays can either be represented 1)~explicitly, resulting in delay differential equations (DDEs)~\cite{Hale:1977}, or 2)~implicitly, by modeling the processes causing the delays. In this paper, we consider optimal control of DDEs, which is relevant to, e.g., the control of chemical and biochemical reactors, power plants, mechanical systems, and systems with communication delays.

The solution to an optimal control problem (OCP) is a set of trajectories of the manipulated inputs, and the objective is either to minimize the deviation from a setpoint or maximize the profit of operation. For general nonlinear systems, it is not possible to derive analytical expressions for the solution of an OCP. Instead, the solution can be approximated numerically~\cite{Binder:etal:2001}. This is referred to as numerical optimal control (NOC). Typically, direct (discretize-then-optimize) or indirect (optimize-then-discretize) approaches are used. In direct approaches, the dynamical constraints are discretized in order to transform the infinite-dimensional OCP into a finite-dimensional nonlinear program (NLP)~\cite{Nocedal:Wright:2006}. In indirect approaches, the optimality conditions are derived analytically and subsequently discretized.
Direct methods are sometimes preferred because they can use existing software for numerical simulation and solution of NLPs. They are divided into 1)~single- and multiple-shooting methods~\cite{Bock:Plitt:1984}, which are based on numerical simulation, and 2)~simultaneous approaches~\cite{Biegler:2007}, where the discretized dynamical equations are constraints in the NLP. A key advantage of multiple-shooting and simultaneous approaches is their ability to handle unstable operating points.

There are many analytical results related to optimal control of DDEs. First-order optimality conditions have been derived for both deterministic~\cite{Boccia:Vinter:2016} and stochastic~\cite{Fuhrman:etal:2010, Oksendal:etal:2011} systems, and specialized methods have been proposed for discretizing the optimality conditions~\cite{Hajipour:Jajarmi:2018}. Furthermore, the Hamilton-Jacobi-Bellman (HJB) equations (for determining optimal feedback control laws) have been derived~\cite{Federico:etal:2010, Federico:etal:2011}, and an efficient approach has been proposed for dynamic programming of DDEs~\cite{Chen:etal:2000}.
Direct approaches have also been considered, and both multiple-shooting~\cite{Bottasso:etal:2009, Scorcelletti:etal:2007} and simultaneous~\cite{Gollmann:etal:2008, Betts:etal:2012} approaches have been proposed.
However, NOC of DDEs is more challenging than for, e.g., ordinary differential equations (ODEs). For instance, if a conventional one-step discretization method is used, either the time step size must be an integer fractions of the delays or the delayed state must be approximated. In both cases, it is difficult to handle time-varying delays, which may depend on both the manipulated inputs and the states. Another challenge is that, in general, the solution to DDEs is not differentiable~\cite{Baker:Paul:1997}. Consequently, the objective function in the OCP may also not be differentiable, and gradient-based methods are not directly applicable to the corresponding NLP.
Therefore, some authors propose to 1)~reformulate the DDEs, e.g., as partial differential equations (PDEs)~\cite{Elnagar:Kazemi:2001}, or 2)~approximate the delayed state, e.g., using lag or Pad\'{e} approximations~\cite{Jiang:etal:2022} or by linearizing the delayed state around the current time~\cite{Gao:etal:2022}. To the best of the authors' knowledge, linearization of the delayed state has not previously been used for NOC of DDEs.

In this work, we present an approach for NOC of DDEs based on linearization of the delayed state variables around the current time. The linearization results in a set of implicit differential equations, and we show that the steady state is identical to that of the original system. However, the stability criteria are different. Therefore, as the approximate system may be unstable, we use a simultaneous approach, and we discretize the differential equations using Euler's implicit method. Finally, we demonstrate the results of the stability analysis and the efficacy of the proposed approach using a numerical example involving a molten salt nuclear fission reactor~\cite{Wooten:Powers:2018, Duderstadt:Hamilton:1976}.

The remainder of the paper is structured as follows. In Section~\ref{sec:optimal:control:problem}, we present the OCP considered in this paper, and in Section~\ref{sec:linearization}, we discuss the linearization of the delayed state. In Section~\ref{sec:simultaneous:approach}, we describe the simultaneous approach, and we present the model of the nuclear reactor in Section~\ref{sec:nuclear:fission}. Finally, the numerical example is presented in Section~\ref{sec:numerical:example}, and conclusions are given in Section~\ref{sec:conclusions}.

	\section{OPTIMAL CONTROL PROBLEM}\label{sec:optimal:control:problem}
We consider optimal control problems in the form
\begin{subequations}\label{eq:ocp}
	\begin{align}\label{eq:ocp:obj}
		\min_{\{u_k\}_{k=0}^{N-1}} \quad \phi &= \phi_x + \phi_{\Delta u}
	\end{align}
	subject to
	\begin{align}
		\label{eq:ocp:ic}
		x(t) &= x_0(t),                            \quad t \in [t_0 - \tau_{\max}, t_0], \\
		\label{eq:ocp:dde}
		\dot x(t) &= f(x(t), z(t), u(t), d(t), p), \quad t \in [t_0, t_f], \\
		\label{eq:ocp:z}
		z(t) &=
		\begin{bmatrix}
			r_1(t - \tau_1) \\
			\vdots \\
			r_m(t - \tau_m)
		\end{bmatrix}, \quad t \in [t_0, t_f], \\
		\label{eq:ocp:r}
		r_i(t) &= h_i(x(t), p), \quad i = 1, \ldots, m, \quad t \in [t_0, t_f], \\
		\label{eq:ocp:u}
		u(t) &= u_k,      \quad t \in [t_k, t_{k+1}[, \quad k = 0, \ldots, N-1, \\
		\label{eq:ocp:d}
		d(t) &= d_k, \quad t \in [t_k, t_{k+1}[, \quad k = 0, \ldots, N-1, \\
		\label{eq:ocp:bounds:x}
		x_{\min} &\leq x(t) \leq x_{\max}, \quad t \in [t_0, t_f], \\
		\label{eq:ocp:bounds:u}
		u_{\min} &\leq u_k \leq u_{\max}, \quad k = 0, \ldots, N-1,
	\end{align}
\end{subequations}
where the objective function, $\phi$, is the sum of a Lagrange term and an input rate-of-movement regularization term, given by
\begin{subequations}
	\begin{align}
		\label{eq:ocp:obj:x}
		\phi_x &= \int_{t_0}^{t_f} \Phi(x(t), u(t), d(t), p) \incr t, \\
		\label{eq:ocp:obj:u}
		\phi_{\Delta u} &= \frac{1}{2} \sum_{k=0}^{N-1} \frac{\Delta u_k^T W_k \Delta u_k}{\Delta t}.
	\end{align}
\end{subequations}
Here, $t$ denotes time, $t_0$ and $t_f$ are the initial and final time, and there are $N$ equidistant control intervals of size $\Delta t$ with the boundaries $t_0 < t_1 < \cdots < t_N = t_f$. Furthermore, $x$ are the states, $z$ are the memory states, $r_i$ for $i = 1, \ldots, m$ are the delayed quantities, $u$ are the manipulated inputs, $d$ are the disturbance variables, and $p$ are the parameters. The constraint~\eqref{eq:ocp:ic} is an initial condition where $x_0$ is the initial state function, \eqref{eq:ocp:dde} are the dynamical constraints (i.e., the DDEs), \eqref{eq:ocp:z}--\eqref{eq:ocp:r} is the contribution from the past states, \eqref{eq:ocp:u}--\eqref{eq:ocp:d} are zero-order-hold (ZOH) parametrizations of the manipulated inputs and disturbance variables, and~\eqref{eq:ocp:bounds:x}--\eqref{eq:ocp:bounds:u} are bounds on the states and manipulated inputs. The DDEs involve $m$ time delays which are functions of the manipulated inputs, i.e., $\tau_i = \tau_i(u(t))$ for $i = 1, \ldots, m$. Consequently, in~\eqref{eq:ocp:ic}, we assume that the initial state function is known for $t \in [t_0 - \tau_{\max}, t_0]$, where
\begin{align}
	\tau_{\max} = \max_{i \in \{1, \ldots, m\}} \max_{u_{\min} \leq u \leq u_{\max}} \tau_i(u).
\end{align}
Depending on the stage cost, $\Phi$, the term in~\eqref{eq:ocp:obj:x} either represents the deviation from a setpoint or the economy of operation. Finally, $W_k$ is a symmetric positive definite weight matrix, and the regularization term in~\eqref{eq:ocp:obj:u} penalizes the change in the manipulated inputs,
\begin{align}
	\Delta u_k &= u_k - u_{k-1}.
\end{align}
The first term in the regularization term~\eqref{eq:ocp:obj:u} involves a reference input, $u_{-1}$, which is given.

\subsection{Steady state and stability}
The solution to OCPs will often drive the system to a steady state, and it is important whether it is stable or not.
In steady state, $\dot x = 0$ and $x(t - \tau_i) = x(t)$, and for the DDEs~\eqref{eq:ocp:dde}--\eqref{eq:ocp:r}, the steady state satisfies
\begin{subequations}\label{eq:dde:steady:state}
	\begin{align}
		0 &= f(x_s, z_s, u_s, d_s, p), \\
		z_s &=
		\begin{bmatrix}
			r_{1, s} \\
			\vdots \\
			r_{m, s}
		\end{bmatrix}, \\
		r_{i, s} &= h_i(x_s, p), \quad i = 1, \ldots, m.
	\end{align}
\end{subequations}
The steady state is locally asymptotically stable if all roots of the characteristic function have strictly negative real part~\cite[Chap.~1, Theorem~6.2]{Hale:1977}. Each root, $\lambda$, is a solution to the characteristic equation
\begin{align}\label{eq:dde:stability}
	\det\left(\lambda I - \pdiff{f}{x} - \sum_{i=1}^m \pdiff{f}{z} \pdiff{z}{r_i} \pdiff{r_i}{x} e^{-\tau_i(u_s) \lambda}\right) &= 0,
\end{align}
where the Jacobian matrices are evaluated in the steady state. In general, there are infinitely many solutions to this equation. We refer to, e.g.,~\cite{Engelborghs:etal:2000} for more details.

	\section{DELAY LINEARIZATION}\label{sec:linearization}
We linearize the delayed states around the current time, $t$:
\begin{align}\label{eq:delay:approximation}
	x(t - \tau_i) &\approx x(t) - \dot x(t) \tau_i.
\end{align}
The resulting approximate system is
\begin{subequations}\label{eq:ide}
	\begin{align}
		\label{eq:ide:x}
		\dot x(t) &= f(x(t), z(t), u(t), d(t), p), \\
		\label{eq:ide:z}
		z(t) &=
		\begin{bmatrix}
			r_1(t) \\
			\vdots \\
			r_m(t)
		\end{bmatrix}, \\
		\label{eq:ide:r}
		r_i(t) &= h(x(t) - \dot x(t) \tau_i, p), \quad i = 1, \ldots, m.
	\end{align}
\end{subequations}
The approximation of the delayed state does not contain any time delays, but it depends on the time derivative of the states. Consequently, the approximate system~\eqref{eq:ide} is a set of implicit differential equations.

\subsection{Steady state and stability}
In steady state, the approximation~\eqref{eq:delay:approximation} is exact because $\dot x = 0$ and $x(t -\tau_i) = x(t)$, and the steady state of~\eqref{eq:ide} also satisfies~\eqref{eq:dde:steady:state}, i.e., it is the same as for the original system.
The steady state is locally asymptotically stable if all roots of the characteristic function have negative real parts~\cite[Prop.~2.1]{Du:etal:2013}. However, the characteristic function is not the same as for the original system. Specifically, each root, $\lambda$, satisfies
\begin{align}\label{eq:ide:stability}
	\det\left(\lambda I - \pdiff{f}{x} - \sum_{i=1}^m \pdiff{f}{z} \pdiff{z}{r_i} \pdiff{r_i}{x} (1 - \tau_i(u_s) \lambda)\right) &= 0,
\end{align}
where the Jacobians are evaluated in the steady state. Note that~\eqref{eq:ide:stability} is also obtained by linearizing the exponential function in~\eqref{eq:dde:stability}, i.e., by using the approximation $e^{-x} \approx 1 - x$. Finally, the characteristic equation~\eqref{eq:ide:stability} corresponds to a generalized eigenvalue problem~\cite[Sect.~7.7]{Golub:VanLoan:2013}.

	\section{SIMULTANEOUS APPROACH}\label{sec:simultaneous:approach}
In this section, we present the simultaneous approach for approximating the solution to the OCP~\eqref{eq:ocp} based on linearization of the delayed state, as described in Section~\ref{sec:linearization}. We discretize the approximate system~\eqref{eq:ide} and the Lagrange term in the objective function~\eqref{eq:ocp:obj:x}, and we present the resulting NLP, which is solved numerically using off-the-shelf software. Furthermore, we present analytical expressions for the first-order derivatives, and we briefly discuss implementation details.

\subsection{Discretization}
We use Euler's implicit method with $M$ time steps per control interval to discretize the differential equations, and we write the discretized equations in residual form. Specifically, the residual equations
\begin{align}\label{eq:residual:function}
	R_{k, n}
	&= R_{k, n}(x_{k, n+1}, x_{k, n}, u_k, d_k, p) \nonumber \\
	&= x_{k, n+1} - x_{k, n} - f(x_{k, n+1}, z_{k, n+1}, u_k, d_k, p) \Delta t_{k, n} \nonumber \\
	&= 0,
\end{align}
must be satisfied for $n = 0, \ldots, M-1$ and $k = 0, \ldots, N-1$. Here, $\Delta t_{k, n}$ is the $n$'th time step size in the $k$'th control interval. Next, we introduce the auxiliary variable
\begin{align}
	\label{eq:discretized:variables:v}
	v_{i, k, n+1} &= x_{k, n+1} - \frac{x_{k, n+1} - x_{k, n}}{\Delta t_{k, n}} \tau_i(u_k),
\end{align}
which approximates the linearized delayed state. Specifically, we use a backward difference approximation of the time derivative in~\eqref{eq:delay:approximation}.
The discretized memory states and delayed quantities are
\begin{subequations}\label{eq:discretized:variables}
	\begin{align}
		\label{eq:discretized:variables:z}
		z_{k, n+1} &=
		\begin{bmatrix}
			r_{1, k, n+1} \\
			\vdots \\
			r_{m, k, n+1}
		\end{bmatrix}, \\
		\label{eq:discretized:variables:r}
		r_{i, k, n+1} &= h_i(v_{i, k, n+1}, p),
	\end{align}
\end{subequations}
for $i = 1, \ldots, m$.
Furthermore, we implicitly enforce the continuity conditions
\begin{subequations}\label{eq:continuity:constraints}
	\begin{align}
		x_{0, 0} &= x_0, \\
		x_{k, 0} &= x_{k-1, M}, & k &= 1, \ldots, N-1,
	\end{align}
\end{subequations}
i.e., we use them to eliminate $x_{k, 0}$ for $k = 0, \ldots, N-1$.
Finally, we use a right rectangle rule to approximate the integral in the Lagrange term~\eqref{eq:ocp:obj:x}:
\begin{align}\label{eq:obj:discretized}
	\psi_x &= \sum_{k=0}^{N-1} \sum_{n=0}^{M-1} \Phi(x_{k, n+1}, u_k, d_k, p) \Delta t_{k, n}.
\end{align}

\subsection{Nonlinear program}
In the simultaneous approach, the solution to the OCP~\eqref{eq:ocp} is approximated by the solution to the NLP
\begin{subequations}\label{eq:nlp}
	\begin{align}\label{eq:nlp:obj}
		\min_{\{\{x_{k, n+1}\}_{n=0}^{M-1}, u_k\}_{k=0}^{N-1}} \quad & \psi = \psi_x + \phi_{\Delta u}
	\end{align}
	subject to
	\begin{align}
		\label{eq:nlp:residual:equations}
		&R_{k, n}(x_{k, n+1}, x_{k, n}, u_k, d_k, p) = 0, \\
		\label{eq:nlp:bounds:x}
		&x_{\min} \leq x_{k, n+1} \leq x_{\max}, \\
		\label{eq:nlp:bounds:u}
		&u_{\min} \leq u_k \leq u_{\max},
	\end{align}
\end{subequations}
where $k = 0, \ldots, N-1$ and $n = 0, \ldots, M-1$. The decision variables are the discretized states and the manipulated inputs. The objective function, $\psi$, in~\eqref{eq:nlp:obj} is the sum of the discretized Lagrange term~\eqref{eq:obj:discretized} and the original regularization term~\eqref{eq:ocp:obj:u}. The constraints~\eqref{eq:nlp:residual:equations} are the discretized differential equations where the continuity constraints~\eqref{eq:continuity:constraints} are enforced implicitly, and in~\eqref{eq:nlp:bounds:x}, we enforce the bounds on the state variables in the discretization points. Finally,~\eqref{eq:nlp:bounds:u} are the original bounds on the manipulated inputs.

\subsection{Jacobian of residual functions}
The Jacobians of the residual functions in~\eqref{eq:residual:function} with respect to the decision variables in the NLP~\eqref{eq:nlp} are
\begin{subequations}
	\begin{align}
		\pdiff{R_{k, n}}{x_{k, n+1}} &=	I - \left(\pdiff{f}{x} + \pdiff{f}{z} \pdiff{z_{k, n+1}}{x_{k, n+1}}\right) \Delta t_{k, n}, \\
		\pdiff{R_{k, n}}{x_{k, n}} &= -I - \pdiff{f}{z} \pdiff{z_{k, n+1}}{x_{k, n}} \Delta t_{k, n}, \\
		\pdiff{R_{k, n}}{u_k} &= -\left(\pdiff{f}{u} + \pdiff{f}{z} \pdiff{z_{k, n+1}}{u_k}\right) \Delta t_{k, n}.
	\end{align}
\end{subequations}
For brevity of the presentation, we omit the arguments of the Jacobians.
The Jacobians of the memory states, $z_{k, n+1}$, and the delayed quantities, $r_{i, k, n+1}$, are in the form
\begin{align}
	\pdiff{z_{k, n+1}}{w} &=
	\begin{bmatrix}
		\pdiff{r_{1, k, n+1}}{w} \\
		\vdots \\
		\pdiff{r_{m, k, n+1}}{w}
	\end{bmatrix}, &
	\pdiff{r_{i, k, n+1}}{w} &= \pdiff{h_i}{x} \pdiff{v_{i, k, n+1}}{w},
\end{align}
where $w$ is $x_{k, n+1}$, $x_{k, n}$, or $u_k$.
The Jacobians of the auxiliary variable in~\eqref{eq:discretized:variables:v} are
\begin{subequations}
	\begin{align}
		\pdiff{v_{i, k, n+1}}{x_{k, n+1}} &= \left(1 - \frac{\tau_i(u_k)}{\Delta t_{k, n}}\right) I, \\
		\pdiff{v_{i, k, n+1}}{x_{k, n}}   &= \frac{\tau_i(u_k)}{\Delta t_{k, n}} I, \\
		\pdiff{v_{i, k, n+1}}{u_k}        &= -\frac{x_{k, n+1} - x_{k, n}}{\Delta t_{k, n}} \pdiff{\tau_i}{u}.
	\end{align}
\end{subequations}

\subsection{Gradient of the objective function}
The gradients of the objective function in~\eqref{eq:nlp:obj} are
\begin{subequations}
	\begin{align}
		\nabla_{x_{k, n+1}} \psi &= \nabla_{x_{k, n+1}} \psi_x, \\
		\nabla_{u_k} \psi &= \nabla_{u_k} \psi_x + \nabla_{u_k} \phi_{\Delta u},
	\end{align}
\end{subequations}
where the gradients of the discretized Lagrange term are
\begin{subequations}
	\begin{align}
		\nabla_{x_{k, n+1}} \psi_x &= \nabla_x \Phi (x_{k, n+1}, u_k, d_k, p) \Delta t_{k, n}, \\
		\nabla_{u_k} \psi_x &= \sum_{n=0}^{M-1} \nabla_u \Phi(x_{k, n+1}, u_k, d_k, p) \Delta t_{k, n},
	\end{align}
\end{subequations}
and the gradients of the regularization term are
\begin{subequations}
	\begin{align}
		\nabla_{u_k} \phi_{\Delta u} &= \frac{W_k \Delta u_k - W_{k+1} \Delta u_{k+1}}{\Delta t}, \\
		\nabla_{u_{N-1}} \phi_{\Delta u} &= \frac{W_{N-1} \Delta u_{N-1}}{\Delta t},
	\end{align}
\end{subequations}
for $k = 0, \ldots, N-2$.

\subsection{Implementation}
We use the interior point algorithm implemented in Matlab's \texttt{fmincon} to solve the NLP~\eqref{eq:nlp}, and we provide the gradient of the objective function and the Jacobian associated with the nonlinear equality constraints (i.e., the Jacobian of the residual functions). We implement the Jacobian as a sparse matrix, and we use a finite-difference approach to approximate the Hessian.

	\section{NUCLEAR FISSION}\label{sec:nuclear:fission}
In this section, we present a model of a molten salt nuclear fission reactor where the salt is circulated through a heat exchanger outside of the reactor core. The model consists of a set of mass and energy balance equations.

The mass balance equations describe the temporal evolution of the concentrations of $N_g = 6$ neutron precursor groups ($i = 1, \ldots, N_g$) and of the neutrons in the reactor ($i = n = N_g+1$):
\begin{subequations}\label{eq:nuclear:fission:concentrations}
	\begin{align}
		\dot C_i(t) &= (C_{i, in}(t) - C_i(t)) D(t) + R_i(t), & i &= 1, \ldots, N_g, \\
		\dot C_n(t) &= R_n(t).
	\end{align}
\end{subequations}
The dilution rate and the volumetric inlet and outlet flow rate are
\begin{align}
	D(t) &= F(t)/V, & F(t) = A v(t),
\end{align}
where $V$ is the reactor volume, $A$ is the cross-sectional area of the inlet and outlet pipe, and $v$ is the velocity of the circulated molten salt. The neutron precursors decay while the molten salt is circulated outside the core. Therefore, the inlet concentration is
\begin{align}\label{eq:nuclear:fission:inlet:concentrations}
	C_{i, in}(t) &= C_i(t - \tau(t)) e^{-\lambda_i \tau(t)},
\end{align}
where $\lambda_i$ is the decay rate of the $i$'th neutron precursor group, and $\tau$ is the time delay, which depends on the length, $L$, of the pipe outside of the reactor core:
\begin{align}
	\tau(t) &= \frac{L}{v(t)}.
\end{align}
The production rates of the neutrons and the neutron precursor groups are
\begin{align}
	R(t) &= S^T\!(t) r(t),
\end{align}
where the stochiometric matrix and the reaction rates are
{\small
\renewcommand{\arraystretch}{1.2}
\begin{align}
	S(t) &=
	\begin{bmatrix}
		-1 & & & 1 \\
		& \ddots & & \vdots \\
		& & -1 & 1 \\
		\beta_1 & \cdots & \beta_{N_g} & \rho(t) - \beta
	\end{bmatrix}\hspace{-2pt}, \;
	r(t) =
	\begin{bmatrix}
		\lambda_1 C_1(t) \\
		\vdots \\
		\lambda_{N_g} C_{N_g}(t) \\
		C_n(t)/\Lambda
	\end{bmatrix}.
\end{align}
}%
Here, $\beta_i$ is the delayed neutron fraction of precursor group~$i$, and $\beta = \sum_{i=1}^{N_g} \beta_i$. Furthermore, $\Lambda$ is the mean neutron generation time, and $\rho$ is the reactivity, which is the sum of the thermal reactivity, $\rho_{th}$, and the external reactivity, $\rho_{ext}$:
\begin{align}
	\rho(t) &= \rho_{th}(t) + \rho_{ext}(t).
\end{align}
The change in the thermal reactivity is negatively proportional to the change in the reactor core temperature, $T_r$:
\begin{align}\label{eq:nuclear:fission:thermal:reactivity}
	\dot \rho_{th}(t) &= -\kappa \dot T_r(t).
\end{align}
The proportionality constant, $\kappa$, is positive.

The energy balances determine the temperature of the reactor core, $T_r$, and the heat exchanger, $T_{hx}$. For both energy balances, we assume that the involved masses, $m_r$ and $m_{hx}$, and the specific heat capacity, $c_P$, are constant in time, and we isolate the temperatures. For the reactor core,
\begin{align}\label{eq:nuclear:fission:reactor:temperature}
	\dot T_r(t) &= \frac{f_r(t)}{m_r} \big(T_{hx}(t - \tau(t)/2) - T_r(t)\big) + \frac{Q_g(t)}{m_r c_P},
\end{align}
where the first term is related to the difference in the enthalpies of the inlet and outlet streams of the core, and the second term is related to the amount of thermal energy generated by nuclear fission in the core, which is proportional to the neutron concentration:
\begin{align}
	Q_g(t) &= Q_{g, 0} \frac{C_n(t)}{C_{n, 0}}.
\end{align}
The nominal thermal energy generation, $Q_{g, 0}$, and neutron concentration, $C_{n, 0}$, are given. The mass flow rate through the core, $f_r$, is equal to the mass flow rate through the heat exchanger, $f_{hx}$:
\begin{align}
	f_r(t) &= f_{hx}(t) = \rho_s A v(t).
\end{align}
Here, $\rho_s$ is the density of the molten salt.
For the heat exchanger,
\begin{align}\label{eq:nuclear:fission:heat:exchanger:temperature}
	\dot T_{hx}(t) &= \frac{f_{hx}(t)}{m_{hx}} \big(T_r(t - \tau(t)/2) - T_{hx}(t)\big) \nonumber \\
	&- \frac{k_{hx}}{m_{hx} c_P} \big(T_{hx}(t) - T_c\big),
\end{align}
where the first term is related to the difference in the enthalpies of the inlet and outlet streams of the heat exchanger, and the second term is related to the amount of energy transferred to the coolant. The conductivity, $k_{hx}$, and the coolant temperature, $T_c$, are constant.

The DDEs described in this section are in the general form~\eqref{eq:ocp:dde}--\eqref{eq:ocp:r}. The states are the $N_g+1$ concentrations given by~\eqref{eq:nuclear:fission:concentrations}, the thermal reactivity in~\eqref{eq:nuclear:fission:thermal:reactivity}, and the reactor and heat exchanger temperatures in~\eqref{eq:nuclear:fission:reactor:temperature} and~\eqref{eq:nuclear:fission:heat:exchanger:temperature}. The manipulated inputs are the velocity, $v$, and the external reactivity, $\rho_{ext}$, there are no disturbance variables, and there are two time delays: One in~\eqref{eq:nuclear:fission:inlet:concentrations} and another in~\eqref{eq:nuclear:fission:reactor:temperature} and~\eqref{eq:nuclear:fission:heat:exchanger:temperature}.

	\section{NUMERICAL EXAMPLE}\label{sec:numerical:example}
\begin{figure*}
	\centering
	\includegraphics[width=0.33\textwidth]{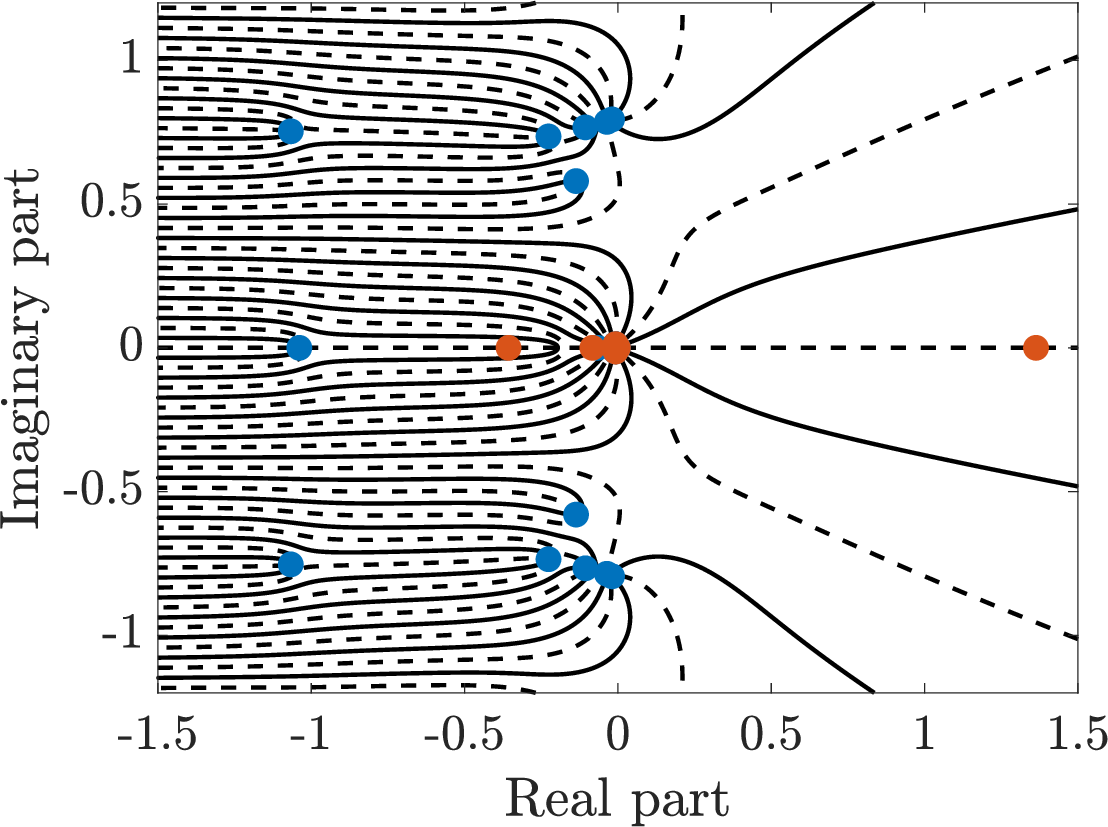}~
	\includegraphics[width=0.33\textwidth]{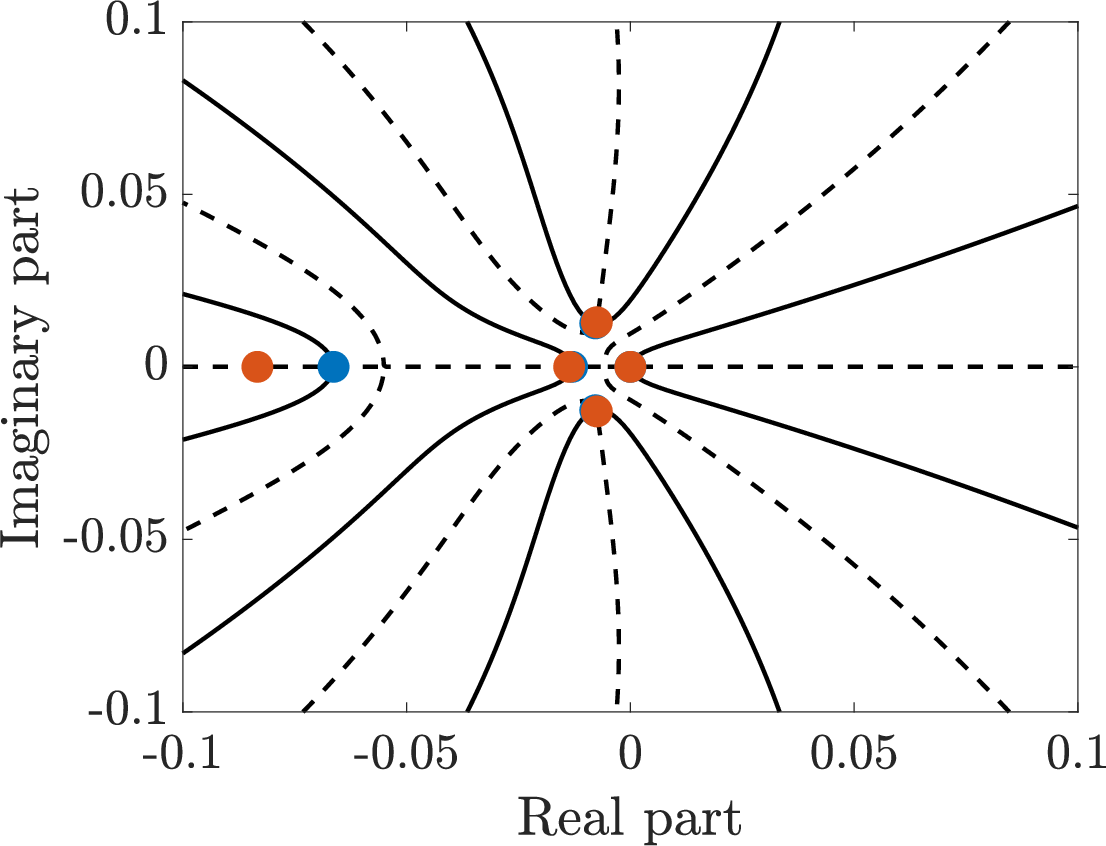}~
	\includegraphics[width=0.3215\textwidth]{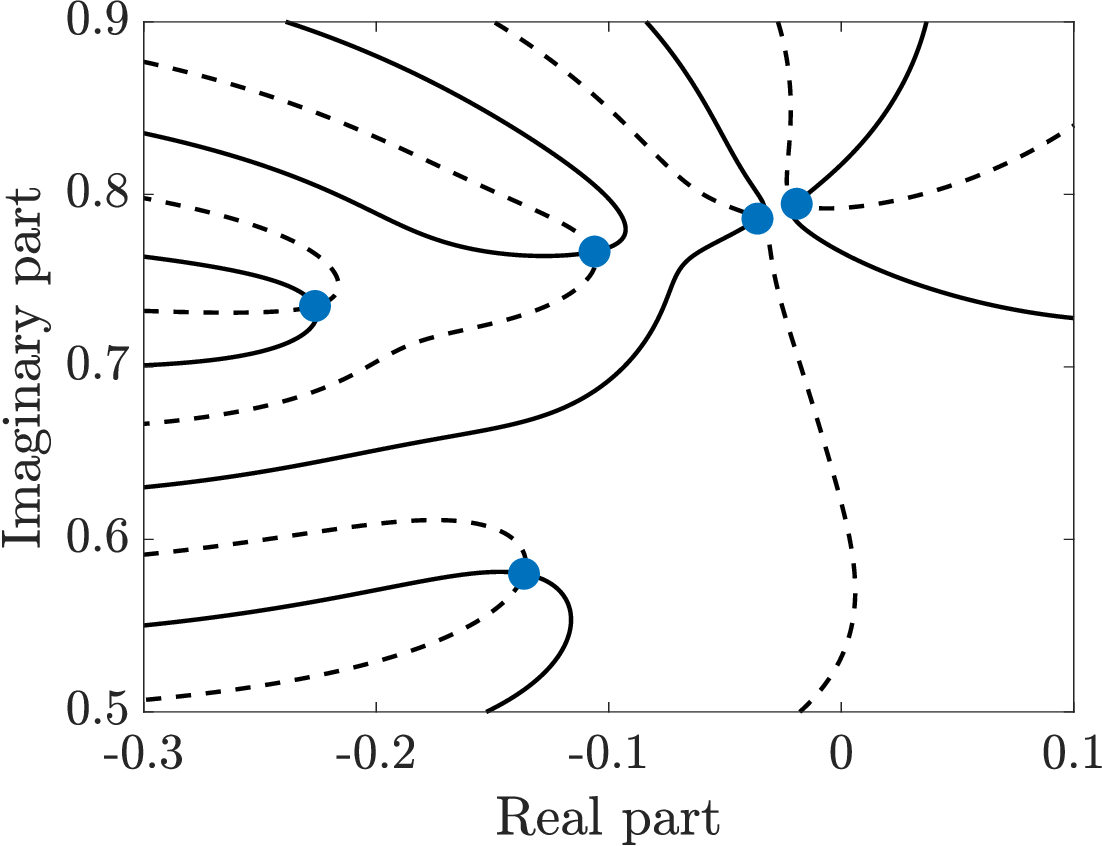}
	\caption{The roots of the characteristic function in~\eqref{eq:dde:stability} for the original DDEs (blue circles), and the roots of the characteristic function in~\eqref{eq:ide:stability} for the approximate system (red circles). The roots $-2.33$, $-4.80$, and $-20.2$ for the approximate system are not shown. The dashed and solid lines indicate the roots of the real and imaginary parts of the characteristic function in~\eqref{eq:dde:stability}, respectively, and their intersections indicate the actual roots.}
	\label{fig:stability:analysis}
\end{figure*}

In this section, we present a numerical example involving the model from Section~\ref{sec:nuclear:fission} with the parameter values listed in Table~\ref{tab:parameters}.
Fig.~\ref{fig:stability:analysis} shows the roots of the characteristic function in~\eqref{eq:dde:stability} for the original DDEs, and the roots of the characteristic function in~\eqref{eq:ide:stability} for the approximate system. Here, $Q_g = 1$~MW, $C_{n, 0} = 1$~m$^{-3}$ (for simplicity), $v = 4$~m/s, and $\rho_{ext} = 50$~pcm (percent mille).
The four roots close to the origin are well approximated by those of the approximate system, but those further away (e.g., those with large imaginary parts) are not.
Furthermore, the approximate system has a positive real root, but we have not identified any roots with positive real part for the original DDEs.

\begin{figure*}[tbh]
	\centering
	\includegraphics[width=\textwidth]{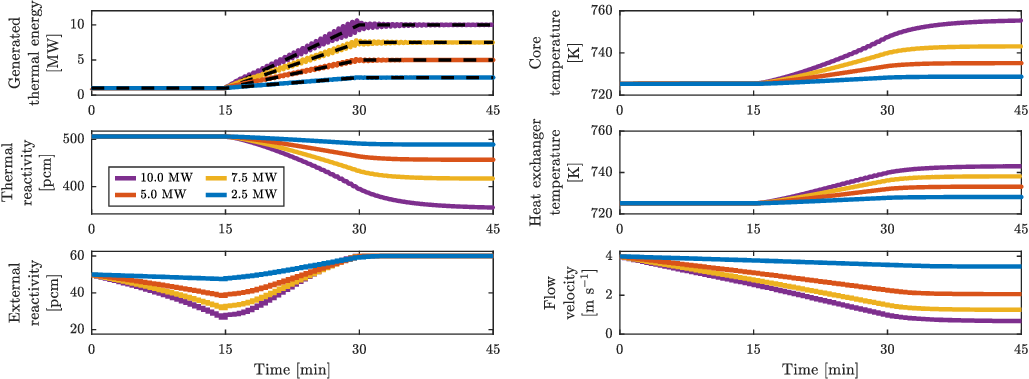}
	\caption{Optimal tracking of four different time-varying setpoints for the thermal energy generation using the simultaneous approach described in Section~\ref{sec:simultaneous:approach}. The top and middle rows show the generated energy, the core and heat exchanger temperatures, and the thermal reactivity when the optimal external reactivity and flow velocity (bottom row) are used to simulate the original system of DDEs presented in Section~\ref{sec:nuclear:fission}.}
	\label{fig:optimal:control:comparison}
\end{figure*}

\begin{figure}[tbh]
	\centering
	\includegraphics[width=0.986\linewidth]{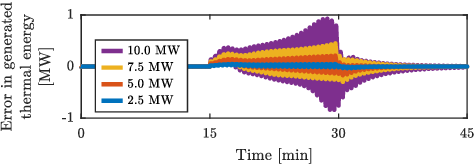}
	\caption{The difference between the generated thermal energy obtained with the original DDEs presented in Section~\ref{sec:nuclear:fission} and with the approximate system.}
	\label{fig:optimal:control:error}
\end{figure}

\begin{figure}[tbh]
	\centering
	\includegraphics[width=0.986\linewidth]{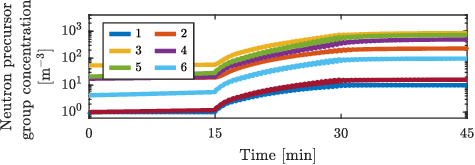}
	\caption{The concentrations of the neutron precursor groups corresponding to the simulation shown in Fig.~\ref{fig:optimal:control:comparison} where the setpoint increases to $10$~MW.}
	\label{fig:optimal:control:neutron:precursor:groups}
\end{figure}

Next, we use the proposed simultaneous approach to track four time-varying setpoints that increase the thermal energy generation from 1~MW to 2.5~MW, 5~MW, 7.5~MW, and 10~MW. The corresponding computation times (for solving the NLP) are 88.2~s, 66.9~s, 124.8~s, and 242.5~s, respectively. The initial guesses of the manipulated inputs are $v = 4$~m/s and $\rho_{ext} = 50$~pcm (this is also the reference input, $u_{-1}$), and for each point in time, the initial guess of the states is the steady state corresponding to the setpoint for the thermal energy generation.
The approximate system is used in the simultaneous approach to compute the manipulated inputs. However, we test the efficacy of the approach by simulating the original DDEs using Matlab's \texttt{ddesd}. The size of the control intervals is $\Delta t = 30$~s, we use $M = 1$ time step per control interval, and $W_k$ is a diagonal matrix with the elements $10^{-2}$~s~pcm$^{-2}$ and $10^2$~s$^3$~m$^{-2}$.
Fig.~\ref{fig:optimal:control:comparison} shows the results, and the reactor is able to track all four setpoints. Furthermore, as the generated power increases, the flow velocity decreases, and the time delays increase. Consequently, the linearization of the delayed states becomes less accurate.
Fig.~\ref{fig:optimal:control:error} shows the difference between the generated thermal energy for the original DDEs and the approximate system. The approximate system is not able to capture the fast oscillations in the generated energy. This is consistent with the results of the stability analysis.
For completeness, Fig.~\ref{fig:optimal:control:neutron:precursor:groups} shows the concentrations of the neutron precursor groups for one of the simulations.
\begin{table}
	\centering
	\caption{Model parameter values.}
	\label{tab:parameters}
	\renewcommand{\arraystretch}{1.2}
	\begin{tabular*}{\linewidth}{@{\extracolsep{\fill}} cccccc}
		\hline
		\multicolumn{6}{c}{Decay constants [$s^{-1}$]} \\
		\hline
		$\lambda_1$ & $\lambda_2$ & $\lambda_3$ & $\lambda_4$ & $\lambda_5$ & $\lambda_6$ \\
		0.0124 & 0.0305 & 0.1110 & 0.3010 & 1.1300 & 3.0000 \\
		\hline
		\multicolumn{6}{c}{Delayed neutron fractions [$-$]} \\
		\hline
		$\beta_1$ & $\beta_2$ & $\beta_3$ & $\beta_4$ & $\beta_5$ & $\beta_6$ \\
		0.00021 & 0.00141 & 0.00127 & 0.00255 & 0.00074 & 0.00027 \\
		\hline
	\end{tabular*}
	\begin{tabular*}{\linewidth}{@{\extracolsep{\fill}} cccc}
		\multicolumn{4}{c}{Other model parameters} \\
		\hline
		$\beta$~[$-$] & $\Lambda$ [s] & $c_P$~[MJ~kg$^{-1}$~K$^{-1}$] & $k_{hx}$~[MW~K$^{-1}$]  \\
		0.0065 & $5\cdot 10^{-5}$ & 2$\cdot 10^{-3}$ & 0.5 \\
		\hline
		$\kappa$~[K$^{-1}$] & $\rho_s$~[kg~m$^{-3}$] & $m_r$~[kg] & $m_{hx}$~[kg] \\
		$5\cdot 10^{-5}$ & 2,000 & 10,000 & 2,500 \\
		\hline
		$V$~[m$^3$] & $A$~[m$^2$] & $L$~[m] & $T_c$~[K] \\
		0.5 & 0.3 & 30 & 723.15 \\
		\hline
	\end{tabular*}
\end{table}

	\section{Conclusions}\label{sec:conclusions}
We propose a simultaneous approach for NOC of DDEs based on linearization of the delayed states. The linearization results in a set of implicit differential equations, which we discretize using Euler's implicit method. Furthermore, we compare the stability criteria of the DDEs and the approximate system, which are also related by a linearization. We use Matlab's \texttt{fmincon} to solve the NLP in the simultaneous approach, and we present a numerical example involving a molten salt nuclear fission reactor. We test the proposed approach using the original DDEs, which successfully track a time-varying setpoint.

	\bibliographystyle{IEEEtran}
	\bibliography{./ref/References}
\end{document}